\documentclass[preprint,3p,12pt,pdf]{elsarticle}

\usepackage{lineno,hyperref}
\modulolinenumbers[5]

\usepackage{hyperref}

\usepackage{epstopdf}

\usepackage{tikz}
\usetikzlibrary{arrows}

\usepackage{pst-node}
\usepackage{tikz-cd}

\usepackage[shellescape]{gmp}

\usepackage{mathrsfs}
\usepackage{amssymb}
\usepackage{amsfonts}
\usepackage{latexsym}
\usepackage{mathtools}
\usepackage{xcolor}
\usepackage{wrapfig}
\usepackage{floatflt}

\usepackage{mathtools}
\usepackage{extarrows}

\usepackage{graphicx}

\usepackage{subcaption}

\def\lb{\label}

\newcommand{\er}[1]{\textrm{(\ref{#1})}}

\newtheorem{theorem}{\bf Theorem}[section]

\def\a{\alpha}         
\def\b{\beta}          
         
\def\G{\Gamma}         
\def\d{\delta}

           \def\mJ{{\mathscr J}}

         \def\mS{{\mathscr S}}
\def\t{\tau}

\def\ve{\varepsilon}   \def\vt{\vartheta}    \def\vp{\varphi}    

\def\Z{{\mathbb Z}}    \def\R{{\mathbb R}}       
    \def\N{{\mathbb N}}   

\def\lt{\biggl}                  \def\rt{\biggr}
               \def\wt{\widetilde}

\let\ge\geqslant                 \let\le\leqslant

\def\iy{\infty}

                 \def\ev{\equiv}
        
\def\el2{\ell^{\,2}}             \def\1{1\!\!1}

\def\arg{\mathop{\mathrm{arg}}\nolimits}

\def\Im{\mathop{\mathrm{Im}}\nolimits}

\def\Re{\mathop{\mathrm{Re}}\nolimits}

\def\sign{\mathop{\mathrm{sign}}\nolimits}

\let\ge\geqslant
\let\le\leqslant

\newcommand{\ca}{\begin{cases}}
	\newcommand{\ac}{\end{cases}}
\newcommand{\ma}{\begin{pmatrix}}
	\newcommand{\am}{\end{pmatrix}}
\renewcommand{\[}{\begin{equation}}
	\renewcommand{\]}{\end{equation}}
\def\eq{\begin{equation}}
	\def\qe{\end{equation}}

\begin{document}
	
	\begin{frontmatter}

		\title{Complete left-tail asymptotic for branching processes with immigration}

		\date{\today}

		\author
		{Anton A. Kutsenko}
	
\address{University of Hamburg, MIN Faculty, Department of Mathematics, 20146 Hamburg, Germany; email: akucenko@gmail.com}
	
	\begin{abstract}
		We derive a complete left-tail asymptotic series for the density of the {\it martingale limit} of a Galton-Watson process with immigration. We show that the series converges everywhere, not only for small arguments. This is the first complete result regarding the left tails of branching processes with immigration. A good, quickly computed approximation for the density will also be derived from the series.
	\end{abstract}

	\begin{keyword}
	 Galton-Watson process, left-tail asymptotic, 
        Schr\"oder and Poincar\'e-type functional equations, Karlin-McGregor function, Fourier analysis
	\end{keyword}

	
\end{frontmatter}


{\section{Introduction}\lb{sec0}}

Let us start with the classical branching process without immigration. Analyzing the tail asymptotics of the density of the {\it martingale limit} of supercritical Galton–Watson processes began in \cite{H}. We will discuss the Schr\"oder case, where one individual's survival probability is non-zero. For the small argument of the density, analyzing the first asymptotic term was started in \cite{D1}, where precise estimates were obtained. Then in \cite{BB1}, the authors found an explicit form of the first asymptotic term. After that, in \cite{K24}, the complete left-tail asymptotic was derived. All the asymptotic terms were expressed explicitly through the so-called one-periodic Karlin-McGregor function introduced in \cite{KM1} and \cite{KM2}. Moreover, it was shown that the series converges everywhere, not only for small arguments of the density. For the Galton–Watson processes with immigration, results are available only for the first asymptotic terms; see, e.g., \cite{CLR} and \cite{BGMS}. We will extend the results from \cite{K24} to the case with immigration. It is impossible to apply directly the results of \cite{K24} to this case, and we need to develop new techniques.

The Galton-Watson process with immigration is defined by
\[\lb{001}
 X_{t+1}=\sum_{j=1}^{X_t}\xi_{j,t}+Y_t,\ \ \ X_0=1,\ \ \ t\in\N\cup\{0\},
\]
where all $\xi_{j,t}$ are independent and identically distributed natural number-valued random variables with the probability-generating function
\[\lb{002}
 P(z):=\mathbb{E}z^{\xi}=p_1z+p_2z^2+p_3z^3+....
\]
We consider the so-called Schr\"oder case $p_0=0$ and $1<p_1<1$. For simplicity, we assume that $P(z)$ is regular at $z=1$ - in particular, the expectation $E=P'(1)$ is finite. All immigration terms $Y_t$ are also identically distributed and independent of other terms. Their probability-generating function is
\[\lb{003}
 Q(z):=\mathbb{E}z^{Y}=q_0+q_1z+q_2z^2+....
\]
We assume that $ q_0\ne 0$ - sometimes, immigrants may not arrive. Again, for simplicity, we require that $Q(z)$ is regular at $z=1$.

To formulate the main result, we extract some information from the next Section. For the classical case $Q\ev1$, the density $p(x)$ of the {\it martingale limit} $E^{-t}X_t$ can be computed as a Fourier transform of some special function satisfying a Poincar\'e-type functional equation:
$$
 p(x)=\frac1{2\pi}\int_{-\iy}^{+\iy}\Pi(\mathbf{i}y)e^{\mathbf{i}yx}dy,
$$
where $\Pi(z)$ is given by
$$
\Pi(z):=\lim_{t\to+\iy}\underbrace{P\circ...\circ P}_{t}(1-\frac{z}{E^{t}}),
$$
it satisfies
$$
P(\Pi(z))=\Pi(Ez),\ \ \ \Pi(0)=1,\ \ \ \Pi'(0)=-1.
$$
This functional equation is fundamental in the analysis of $p(x)$. Defining another analytic function
$$
 \Phi(z):=\lim_{t\to+\iy}p_1^{-t}\underbrace{P\circ...\circ P}_{t}(z),
$$
satisfying the Schr\"oder-type functional equation
$$
  \Phi(P(z))=p_1\Phi(z),\ \ \ \Phi(z)\sim z,\ z\to0,
$$
one can use the one-periodic Karlin-McGregor function
$$
 K(z):=p_1^{-z}\Phi(\Pi(E^z))
$$
to obtain the complete left-tail series for $p(x)$, that converges everywhere, not only for small $x>0$. This research was done in \cite{K24}.

For the Galton-Watson process with immigration, the density $p_{\rm imm}(x)$ of the {\it martingale limit} $E^{-t}X_t$ can be computed as a Fourier transform of another special function
\[\lb{T1f4}
  p_{\rm imm}(x)=\frac1{2\pi}\int_{-\iy}^{+\iy}\Pi_{\rm imm}(\mathbf{i}y)e^{\mathbf{i}yx}dy,
\]
where
$$
 \Pi_{\rm imm}(z)=\Pi(z)\cdot R(z),
$$
with
$$
 R(z)=\prod_{t=1}^{+\iy}Q(\Pi(\frac{z}{E^t})) 
$$
satisfying the functional equation
$$
 R(Ez)=R(z)Q(\Pi(z)),\ \ \ R(0)=1.
$$
There are no more Poincar\'e-type functional equations. However, one may still define a one-periodic function
$$
 L(z):=q_0^{-z}{R(E^z)}{\Psi(\Pi(E^z))},
$$
with the help of another analytic function
$$
 \Psi(z):=\prod_{t=0}^{+\iy}q_0^{-1}Q(\underbrace{P\circ...\circ P}_{t}(z)),
$$
satisfying another functional equation
$$
  Q(z)\Psi(P(z))=q_0\Psi(z),\ \ \ \Psi(0)=1.
$$
To formulate our main result, we need to introduce a few more objects. The first one is the {\it critical angle} $\wt\theta^*$ - this is the maximal angle of a sector of arbitrarily small radius in the neighborhood of $1$ that lies entirely in the filled Julia set for the function $P(z)$. Another object is the Fourier coefficients of the following one-periodic functions
$$
 (K^{n+1}\cdot L)(z)=\sum_{m=-\iy}^{m=+\iy}\vt_{nm}e^{2\pi\mathbf{i}mz},\ \ \ n\ge0. 
$$
We also need the Taylor coefficients of the analytic function
\[\lb{T1f5}
 A(z)=\frac{\Phi^{-1}(z)}{z\Psi(\Phi^{-1}(z))}=\sum_{n=0}^{+\iy}A_nz^n.
\]
\begin{theorem}\lb{T1} If $\wt\theta^*>\pi$ and $\log_E(p_1q_0)<-1$ then
\[\lb{T1f1}
 p_{\rm imm}(x)=\sum_{n=0}^{+\iy}A_nx^{-\log_E(p_1^{n+1}q_0)-1}B_n(-\log_Ex),\ \ \ x>0,
\]
where one-periodic functions $B_n(z)$ are defined by
\[\lb{T1f2}
 B_n(z)=\sum_{m=-\iy}^{m=+\iy}\frac{\vt_{nm}e^{2\pi\mathbf{i}mz}}{\G(-\frac{\ln(p_1^{n+1}q_0)+2\pi\mathbf{i}m}{\ln E})}.
\]
Substituting \er{T1f2} into \er{T1f1}, we get a double series that converges super-exponentially fast in $n$ and exponentially fast in $m$, namely we have
\[\lb{T1f3}
 \lt|\frac{A_n\vt_{nm}}{\G(-\frac{\ln(p_1^{n+1}q_0)+2\pi\mathbf{i}m}{\ln E})}\rt|\le C e^{-\a n\ln(n+|m|)-\b|m|},
\]
for some constants $\a,\b,C>0$.
\end{theorem}

{\bf Remark.} The condition $\log_E(p_1q_0)<-1$ can be removed from the Theorem \ref{T1}. The only place where we use it is \er{220}. However, skipping the first finite number of terms, all the other terms give the absolute convergence of the corresponding integrals, since $\log_E(p_1^{n+1}q_0)<-1$ for all large $n$. The accurate explanation of this statement may be left as an interesting exercise for interested readers.

That the {\it critical angle} $\wt\theta^*\ge\pi$ is always true is obvious, since the unit disc belongs to the filled Julia set for any probability-generating function $P(z)$. We demand a little more $\wt\theta^*>\pi$. However, this requirement can also be weakened in some sense. 

Usually, the Fourier coefficients of $B_n(z)$ are small, except for the zeroth one. The reason for the more precise estimates of the form \er{T1f3} is the width of the strip, where $B_n(z)$ are defined. The width is large enough. Hence, the Fourier coefficients of $B_n(z)$ decay quickly. Such a type of phenomenon is known in physics and biology regarding other objects, see, e.g. \cite{DIL,DMZ,CG}.
This fact allows us to define a good, quickly computed approximation of the density
\[\lb{T1f6}
 \wt p_{\rm imm}(x)=\sum_{n=0}^M\frac{(K^{n+1}\cdot L)(0)A_n}{\G(-\frac{\ln(p_1^{n+1}q_0)+2\pi\mathbf{i}m}{\ln E})}x^{-\log_E(p_1^{n+1}q_0)-1},\ \ \ x>0,
\]
with some moderate $M$. For the computation $K(0)$ and $L(0)$, one may use fast convergent methods, similar to those developed in, e.g., \cite{K}. The computation of $A_n$ is generally explicit, because $\Phi^{-1}(z)$ satisfies the Poincar\'e-type functional equation
\[\lb{004}
 \Phi^{-1}(p_1z)=P(\Phi^{-1}(z)),\ \ \ \Phi^{-1}(z)\sim z,\ z\to0,
\]
and the composition $\Psi(\Phi^{-1}(z))$ satisfies
\[\lb{005}
 Q(\Phi^{-1}(z))\Psi(\Phi^{-1}(p_1z))=q_0\Psi(\Phi^{-1}(z)),\ \ \ \Psi(\Phi^{-1}(0))=1.
\]
Substituting Taylor expansions in \er{004} and \er{005}, one may find Taylor coefficients for both functions step by step. Then one may use them to find $A_n$ in \er{T1f5}.

\begin{figure}
	\centering
	\begin{subfigure}[b]{0.75\textwidth}
		\includegraphics[width=\textwidth]{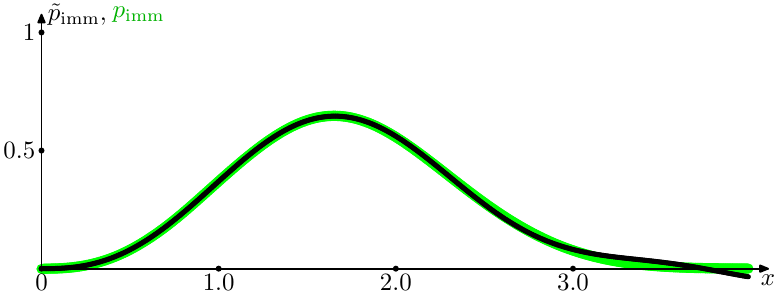}
		\caption{$p_1=0.3$ and $q_0=0.5$}
		\label{fig1a}
	\end{subfigure}
	\hfill
	\begin{subfigure}[b]{0.75\textwidth}
		\includegraphics[width=\textwidth]{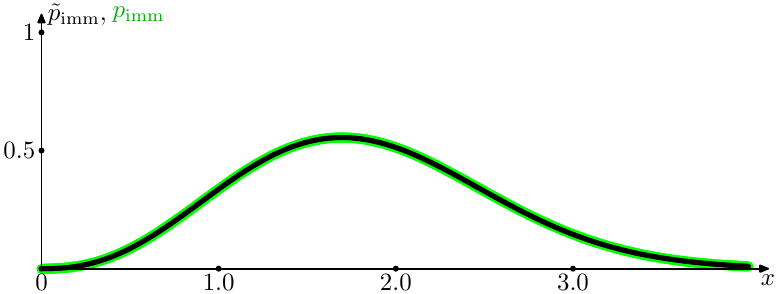}
		\caption{$p_1=0.4$ and $q_0=0.5$}
		\label{fig1b}
	\end{subfigure}
	\hfill
	\begin{subfigure}[b]{0.75\textwidth}
		\includegraphics[width=\textwidth]{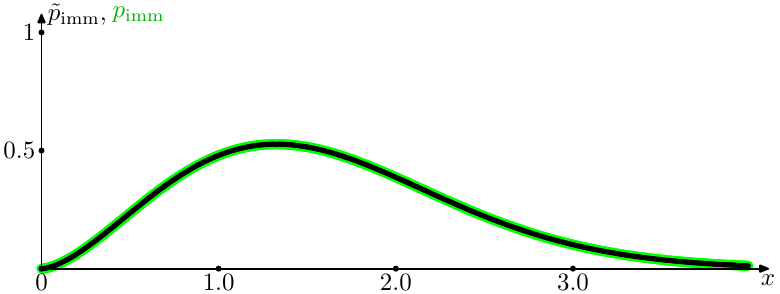}
		\caption{$p_1=0.5$ and $q_0=0.7$}
		\label{fig1c}
	\end{subfigure}
	\caption{For the case \er{006}, the approximation of the density \er{T1f6} is compared with \er{T1f4}.}\label{fig0}
\end{figure}

As an illustration, we compare the density computed by \er{T1f4} and \er{T1f6} for the case
\[\lb{006}
 P(z)=p_1z+(1-p_1)z^2,\ \ \ Q(z)=q_0+(1-q_0)z.
\]
For the computation of \er{T1f4}, we take the interval $y\in[-2000,2000]$ (in fact,  $[0,2000]$, and use the symmetry), devided by $10^6$ points in the trapezoidal rule; the function $\Pi_{\rm imm}$ is computed up to $70$ iterations of $P$. In \er{T1f6}, we take $M=10$. For the calculations, Embarcadero Delphi Rad Studio Community Edition and the library NesLib.Multiprecision are used. This software provides a convenient environment for programming and well-functioning basic functions for high-precision computations. All the algorithms related to the article's subject, including efficient parallelization, are developed by the author (AK). We use 128-bit precision instead of the standard 64-bit double precision. The results of the comparison are quite good, at least for moderate values of $x$, see Fig. \ref{fig0}. In particular, the areas under the plotted curves are $1\pm10^{-2}$. Thus, the density values at which the curves may significantly differ are quite small. The calculation of the approximation \er{T1f6} is real-time, but \er{T1f4} requires a few minutes on Intel Core i7-7700HQ.

{\section{Proof of Theorem \ref{T1}}\lb{sec1}}

\subsection{Galton-Watson process without immigration}

The probability generating function of $X_t$ can be computed as the following composition
\[\lb{100}
 P_t(z)=\underbrace{P\circ...\circ P}_{t}(z)=\sum_{n=1}^{+\iy}\mathbb{P}(X_t=n)z^n.
\]
There are two limit distributions. The first one is not very common but important, the distribution of {\it  ratios of probabilities of rare events}
\[\lb{101}
 \vp_{n}=\lim_{t\to+\iy}\frac{\mathbb{P}(X_t=n)}{\mathbb{P}(X_t=1)}=\lim_{t\to+\iy}p_1^{-t}\mathbb{P}(X_t=n).
\]
The corresponding generating function
\[\lb{102}
 \Phi(z)=\lim_{t\to+\iy}p_1^{-t}P_t(z)=\sum_{n=1}^{+\iy}\vp_{n}z^n
\]
satisfy, by \er{100}-\er{102}, the Schr\"oder type functional equation
\[\lb{103}
 \Phi(P(z))=p_1\Phi(z),\ \ \ \Phi(z)\sim z,\ z\to0.
\]

\begin{figure}[h]
	\center{\includegraphics[width=0.7\linewidth]{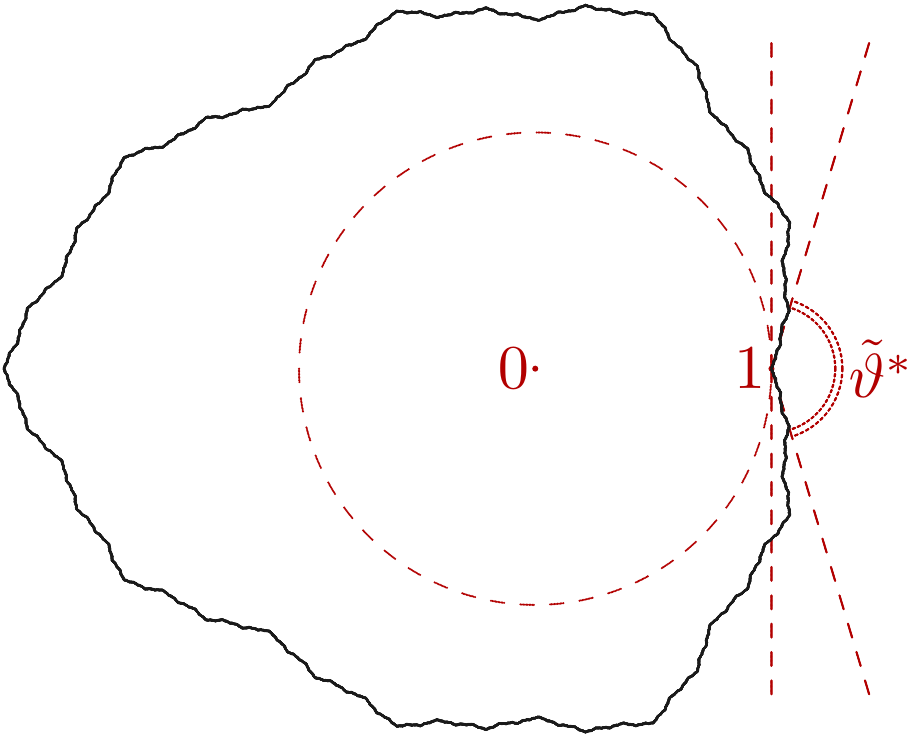}}
	\caption{Julia set for the polynomial $P(z)=0.1z+0.5z^2+0.4z^3$. It is the boundary of the open filled Julia set $\mJ_P$, which is the domain of definition for the analytic function $\Phi$, as seen in Equations \er{102} and \er{103}. The unit disc belongs to the filled Julia set. The figure is taken from my article \cite{K}.}\lb{fig1}
\end{figure}

The function $\Phi(z)$ is analytic inside the open filled Julia set $\mJ_P$ related to $P(z)$, see an example in Fig. \ref{fig1}. This is because $P_t(z)\to0$ when $z\in \mJ_P$, and the convergence to zero is exponentially fast with the factor $p_1=P'(0)$, see details in \cite{K24}.

Another distribution is well-known, the {\it martingale limit} $W=\lim_{t\to+\iy}E^{-t}X_t$. The most important characteristic of $W$ is its density $p(x)$. We follow the ideas of \cite{D1} to obtain the known formula for the density: the Fourier transform of the Poincar\'e function 
\[\lb{104}
 p(x)=\frac1{2\pi}\int_{-\iy}^{+\iy}\Pi(\mathbf{i}y)e^{\mathbf{i}yx}dy,
\]
where $\Pi(z)$ is given by
\[\lb{105}
\Pi(z):=\lim_{t\to+\iy}\underbrace{P\circ...\circ P}_{t}(1-\frac{z}{E^{t}}).
\]
This function is analytic and satisfies the Poincar\'e-type functional equation
\[\lb{106}
P(\Pi(z))=\Pi(Ez),\ \ \ \Pi(0)=1,\ \ \ \Pi'(0)=-1,
\]
which follows from \er{105}. If $P$ is entire, then $\Pi$ is entire as well. 

We give an idea of the derivation of \er{104}, skipping many details available in \cite{D1} and other sources. The density of $W$, roughly speaking, satisfies 
\[\lb{107a}
 p(x)=\lim_{dx\to0}\frac{\mathbb{P}(W\in[x,x+dx))}{dx}
\] 
or
\[\lb{107}
 p(x)=\lim_{t\to+\iy}E^t\mathbb{P}(X_t=[xE^t]),\ x>0,\ \ \ p(x)=0,\ x<0,
\] 
where $[...]$ denotes the integer part. On the one side, from \er{100}, we have the Riemann sum
\[\lb{108}
\underbrace{P\circ...\circ P}_{t}(e^{\frac{-\mathbf{i}y}{E^t}})=E^{-t}\sum_{n=0}^{+\iy} E^{t}\mathbb{P}(X_t=n)e^{\frac{-n\mathbf{i}y}{E^t}}\to \int_{-\iy}^{+\iy}p(x)e^{-\mathbf{i}yx}dx.
\]
On the other side, \er{105} and the well-known expansion for the exponent give
\[\lb{109}
\underbrace{P\circ...\circ P}_{t}(e^{\frac{-\mathbf{i}y}{E^t}})\sim \underbrace{P\circ...\circ P}_{t}(1-{\frac{\mathbf{i}y}{E^t}})\to\Pi(\mathbf{i}y).
\]
Thus, combining \er{108} and \er{109}, we get
\[\lb{110}
 \Pi(\mathbf{i}y)=\int_{-\iy}^{+\iy}p(x)e^{-\mathbf{i}yx}dx.
\] 
Applying the inverse Fourier transform to \er{110}, we arrive at \er{104}.

Due to \er{103} and \er{106}, one may define the one-periodic analytic function
\[\lb{111}
 K(z)=p_1^{-z}\Phi(\Pi(E^z)).
\]
This is the so-called Karlin-McGregor function, which is very helpful in the analysis of $\Pi(z)$ and $p(x)$. From \er{111}, we have 
\[\lb{112}
 \Pi(z)=\Phi^{-1}(z^{\log_Ep_1}K(\log_Ez)).
\]
By definitions \er{105} and \er{111}, the one-periodic function $K(z)$ is defined, at least, in the strip 
\[\lb{113}
 \mS=\{|\Im z|<\frac{\wt\theta^*}{2\ln E}\},
\] 
where $\wt\theta^*$ is the so-called critical angle, see Fig. \ref{fig1} - the maximal angle for which $1-re^{\mathbf{i}\theta}$, with $\theta<\wt\theta^*$, and $r\to+0$ belongs to the filled Julia set. This is because if $z\in\mS$ and $\Re z\to-\iy$ then $\Pi(E^z)=1-E^z+o(E^z)$ belongs to $\mJ_P$, and, hence $K(z)$ in \er{111} is well-defined. Then, one-periodicity of $K(z)$ allows us to extend the domain of definition from $\Re z\to-\iy$ to the whole $\Re z\in\R$, or $z\in\mS$. This explanation is also given in \cite{K24}. In particular, by \er{112}, we have
\[\lb{114}
 \Pi(z)=O(z^{\log_Ep_1}),\ \ \ z\to\iy,\ \arg z<\wt\theta^*-\d,
\]
uniformly for any sufficiently small $\d>0$. Because $\log_Ep_1<0$ and $\Pi(z)$ tends to $0$ for $z\to\iy$ in the corresponding sector, see \er{114}, we can write
\[\lb{115}
p(x)=\frac1{2\pi\mathbf{i}}\int_{\ve-\mathbf{i}\iy}^{\ve+\mathbf{i}\iy}\Pi(z)e^{zx}dz,\ \ \ \ve>0,
\]
by \er{104} and Cauchy's integral theorem in complex analysis - we have shifted the line of integration.

\subsection{Galton-Watson process with immigration}

The difference with the classical Galton-Watson process is the income of new individuals at each time step. Thus, instead of $P_{t+1}(z)=P_t(P(z))$, see \er{100}, we have
\[\lb{200}
 P_{{\rm imm},t+1}(z)=P_{{\rm imm},t}(P(z))Q(z),
\]
which gives
\[\lb{201}
 P_{{\rm imm},t}(z)=\underbrace{P\circ...\circ P}_{t}(z)\cdot Q(\underbrace{P\circ...\circ P}_{t-1}(z))\cdot...\cdot Q(P(z))\cdot Q(z).
\]
The analog of the generating function for the {\it  ratios of probabilities of rare events}, see \er{102}, is
\[\lb{202}
 \Phi_{\rm imm}(z)=\lim_{t\to+\iy}(p_1q_0)^{-t}P_{{\rm imm},t}(z)=\Phi(z)\Psi(z),
\]
where
\[\lb{203}
 \Psi(z)=\lim_{t\to+\iy}q_0^{-t}Q(\underbrace{P\circ...\circ P}_{t-1}(z))\cdot...\cdot Q(P(z))\cdot Q(z)=\prod_{t=0}^{+\iy}q_0^{-1}Q(P_t(z)),
\]
where $P_0(z)=z$. The function $\Psi(z)$ is analytic inside the filled Julia set $\mJ_P$ related to the probability-generating function $P(z)$, since $q_0^{-1}Q(P_t(z))-1\sim q_0^{-1}q_jP_t(z)^j$ with $j>0$ is minimal for which $q_j\ne0$. Then, using the fact that $P_t(z)\sim p_1^tz$ is exponentially small for large $t$ and $z\in\mJ_P$, we see that the infinite product in \er{203} converges to an analytic function. It is seen from \er{203} that $\Psi(z)$ satisfies the functional equation
\[\lb{204}
 Q(z)\Psi(P(z))=q_0\Psi(z),\ \ \ \Psi(0)=1.
\]
The full analog of $\Pi(z)$, see \er{105} and \er{201}, is
\[\lb{205}
 \Pi_{\rm imm}(z)=\lim_{t\to+\iy}P_{{\rm imm},t}(1-\frac{z}{E^t}),
\]
and, hence, the density, see \er{104}, is
\[\lb{206}
 p_{\rm imm}(x)=\frac1{2\pi}\int_{-\iy}^{+\iy}\Pi_{\rm imm}(\mathbf{i}y)e^{\mathbf{i}yx}dy.
\]
Due to \er{105}, \er{201}, and \er{205}, we have
\[\lb{207}
 \Pi_{\rm imm}(z)=\Pi(z)\cdot R(z),
\]
where
\[\lb{208}
 R(z)=\prod_{t=1}^{+\iy}Q(\Pi(\frac{z}{E^t})) 
\]
satisfies the functional equation
\[\lb{209}
 R(Ez)=R(z)Q(\Pi(z)),\ \ \ R(0)=1.
\]
The series \er{208} converges because 
\[\lb{210}
 Q(\Pi(\frac{z}{E^t}))-1\sim Q'(1)(\Pi(\frac{z}{E^t})-1)\sim Q'(1)\frac{z}{E^t}
\] 
for large $t$. Substituting $\Pi(z)$ instead of $z$ into \er{204} and using \er{106}, we obtain
\[\lb{211}
  Q(\Pi(z))\Psi(\Pi(Ez))=q_0\Psi(\Pi(z)).
\]
Eliminating $Q(\Pi(z))$ from \er{209} and \er{211}, we deduce that
\[\lb{212}
 {R(Ez)}{\Psi(\Pi(Ez))}=q_0{R(z)}{\Psi(\Pi(z))},
\]
which means that 
\[\lb{213}
 L(z):=q_0^{-z}{R(E^z)}{\Psi(\Pi(E^z))}
\]
is a one-periodic function. Thus,
\[\lb{214}
 R(z)=z^{\log_Eq_0}\frac{L(\log_Ez)}{\Psi(\Pi(z))},
\]
and
\[\lb{215}
 \Pi_{\rm imm}(z)=z^{\log_Eq_0}\frac{\Pi(z)L(\log_Ez)}{\Psi(\Pi(z))}.
\]
Substituting \er{112} into \er{215} leads to
\[\lb{216}
 \Pi_{\rm imm}(z)=z^{\log_E(p_1q_0)}(K\cdot L)(\log_Ez)A(z^{\log_Ep_1}K(\log_Ez)),
\]
where
\[\lb{217}
 A(z)=\frac{\Phi^{-1}(z)}{z\Psi(\Phi^{-1}(z))}=\sum_{n=0}^{+\iy}A_nz^n
\]
is a function analytic in some neighborhood of $z=0$, with $A_0=1$. Repeating the arguments explained between \er{112} and \er{115}, we say
\[\lb{218}
\Pi_{\rm imm}(z)=O(z^{\log_E(p_1q_0)}),\ \ \ z\to\iy,\ \arg z<\wt\theta^*-\d,
\]
uniformly for any sufficiently small $\d>0$, and
\[\lb{219}
 p_{\rm imm}(x)=\frac1{2\pi\mathbf{i}}\int_{\ve-\mathbf{i}\iy}^{\ve+\mathbf{i}\iy}\Pi_{\rm imm}(z)e^{zx}dz,\ \ \ \ve>0.
\]
The arguments can be repeated, since \er{207} and \er{208} are fulfilled and $|Q(z)|\le1$ for all sufficiently small $z$, i.e., for $|z|\le1$. Also, the one-periodic function $L(z)$, as well as one-periodic $K(z)$, is analytic in the strip $\mS$, see \er{113}. Thus, substituting \er{216} along with \er{217} into \er{219}, we obtain
\[\lb{220}
 p_{\rm imm}(x)=
 \sum_{n=0}^{+\iy}\frac{A_n}{2\pi\mathbf{i}}\int_{\ve-\mathbf{i}\iy}^{\ve+\mathbf{i}\iy}z^{\log_E(p_1^{n+1}q_0)}(K^{n+1}\cdot L)(\log_Ez)e^{zx}dz
\]
for all sufficiently large $\ve>0$ (the choice of $\ve$ depends on the radius of convergence of \er{217}), since the ``periodic'' component $(K^{n+1}\cdot L)(\log_Ez)$ and $e^{zx}$ are bounded on $\ve+\mathbf{i}\R$ and $z^{\log_E(p_1^{n+1}q_0)}$ tends to $0$ quickly for $z\to\iy$. In particular, all the integrals converge absolutely, since  $\log_E(p_1q_0)<-1$.

To continue, we need the following auxiliary statement. Suppose that $T(z)=\sum_{m\in\Z}\t_me^{2\pi\mathbf{i}mz}$ is an one-periodic function continuous in the strip $|\Im z|<s$, for some $s>0$. Then, for any sufficiently small $\ve>0$, we have
\[\lb{221}
 |\t_m|\le e^{2\pi |m|(s-\ve)}\max_{x\in[0,1]}|T(x+\sign(m)\mathbf{i}(s-\ve))|.
\]
This statement follows easily from 
\[\lb{222}
 \t_m=\int_{[0,1]+\sign(m)\mathbf{i}(s-\ve)}T(z)e^{-2\pi\mathbf{i}mz}dz.
\]
Denote the Fourier coefficients of one-periodic functions $(K^{n+1}\cdot L)(z)$ by
\[\lb{223}
 (K^{n+1}\cdot L)(z)=\sum_{m=-\iy}^{m=+\iy}\vt_{nm}e^{2\pi\mathbf{i}mz},\ \ \ n\ge0. 
\]
Remembering that $K$ and $L$ are analytic in the strip $\mS$, see \er{115}, and applying \er{221} to \er{223}, we obtain
\[\lb{224}
 |\vt_{nm}|\le e^{\frac{\pi |m|(\wt\theta^*-\ve)}{\ln E}} C_{\ve}^{n},
\]
for any sufficiently small $\ve>0$ and some constants $C_{\ve}$ depending on the maximal values of $(K^{n+1}\cdot L)(z)$ at the corresponding intervals within $\mS$. The Fourier coefficients $\vt_{nm}$ decay exponentially fast in $m$, so we can write
\begin{multline}\lb{225}
 \frac{1}{2\pi\mathbf{i}}\int_{\ve-\mathbf{i}\iy}^{\ve+\mathbf{i}\iy}z^{\log_E(p_1^{n+1}q_0)}(K^{n+1}\cdot L)(\log_Ez)e^{zx}dz=\\
\frac{1}{2\pi\mathbf{i}}\int_{\ve-\mathbf{i}\iy}^{\ve+\mathbf{i}\iy}\sum_{m=-\iy}^{m=+\iy}\vt_{nm}z^{\frac{\ln(p_1^{n+1}q_0)+2\pi\mathbf{i}m}{\ln E}}e^{zx}dz=\\
\sum_{m=-\iy}^{m=+\iy}\frac{1}{2\pi\mathbf{i}}\int_{\ve-\mathbf{i}\iy}^{\ve+\mathbf{i}\iy}\vt_{nm}z^{\frac{\ln(p_1^{n+1}q_0)+2\pi\mathbf{i}m}{\ln E}}e^{zx}dz=\\
\sum_{m=-\iy}^{m=+\iy}\frac{\vt_{nm}x^{-\frac{\ln(p_1^{n+1}q_0)+2\pi\mathbf{i}m}{\ln E}-1}}{2\pi\mathbf{i}}\int_{\ve x^{-1}-\mathbf{i}\iy}^{\ve x^{-1}+\mathbf{i}\iy}z^{\frac{\ln(p_1^{n+1}q_0)+2\pi\mathbf{i}m}{\ln E}}e^{z}dz=\\
\sum_{m=-\iy}^{m=+\iy}\frac{\vt_{nm}x^{-\frac{\ln(p_1^{n+1}q_0)+2\pi\mathbf{i}m}{\ln E}-1}}{\G(-\frac{\ln(p_1^{n+1}q_0)+2\pi\mathbf{i}m}{\ln E})}=x^{-\frac{\ln(p_1^{n+1}q_0)}{\ln E}-1}B_n(-\log_Ex),
\end{multline}
where, for the latest integral, we use the Hankel representation of the $\Gamma$-function, see Example 12.2.6 on p. 254 in \cite{WW}. The one-periodic functions $B_n(z)$ are defined by
\[\lb{226}
 B_n(z)=\sum_{m=-\iy}^{m=+\iy}\frac{\vt_{nm}e^{2\pi\mathbf{i}mz}}{\G(-\frac{\ln(p_1^{n+1}q_0)+2\pi\mathbf{i}m}{\ln E})}.
\]
As we already discussed, $\vt_{nm}$ decays to $0$ exponentially fast. While the values of the Gamma function in \er{226} increase super-exponentially fast in $n$ for fixed $m$, they are also exponentially small for large $m$ and fixed $n$. More precisely, Stirling's approximation of the Gamma function
\[\lb{227}
 \G(z)=\sqrt{\frac{2\pi}z}\lt(\frac{z}{e}\rt)^z\lt(1+O\lt(\frac1z\rt)\rt)
\]
gives
\[\lb{228}
 \lt|\G\lt(-\frac{\ln(p_1^{n+1}q_0)+2\pi\mathbf{i}m}{\ln E}\rt)\rt|\ge\frac{Be^{-\frac{\ln(p_1^{n+1}q_0)}{\ln E}\ln\frac{\sqrt{\ln^2(p_1^{n+1}q_0)+4\pi^2m^2}}{\ln E}+\frac{\ln(p_1^{n+1}q_0)}{\ln E}-\frac{\pi^2|m|}{\ln E}}}{\sqrt{\ln^2(p_1^{n+1}q_0)+4\pi^2m^2}},
\]
for some constant $B>0$. Combaining \er{228} with \er{224}, and using the condition that $\wt\theta^*-\ve>\pi$ for small $\ve>0$, we obtain
\[\lb{229}
 \lt|\frac{\vt_{nm}}{\G(-\frac{\ln(p_1^{n+1}q_0)+2\pi\mathbf{i}m}{\ln E})}\rt|\le C_1 e^{-\a_1 n\ln(n+|m|)-\b_1|m|},
\]
for some constants $\a_1,\b_1,C_1>0$. Due to Cauchy's estimates, Taylor's coefficients of $A(z)$, see \er{217}, do not grow faster than an exponent $r^n$ for some $r>0$. Thus, if we multiply \er{229} by $A_n$, the super-exponential decay by $n$ does not change
\[\lb{230}
 \lt|\frac{A_n\vt_{nm}}{\G(-\frac{\ln(p_1^{n+1}q_0)+2\pi\mathbf{i}m}{\ln E})}\rt|\le C e^{-\a n\ln(n+|m|)-\b|m|},
\]
for some constants $\a,\b,C>0$. Combining \er{220}, \er{225}, and \er{226}, we obtain
\[\lb{231}
 p_{\rm imm}(x)=\sum_{n=0}^{+\iy}A_nx^{-\frac{\ln(p_1^{n+1}q_0)}{\ln E}-1}B_n(-\log_Ex),\ \ \ x>0.
\]

\section*{Acknowledgements} 
This paper is a contribution to the project S1 of the Collaborative Research Centre TRR 181 "Energy Transfer in Atmosphere and Ocean" funded by the Deutsche Forschungsgemeinschaft (DFG, German Research Foundation) - Projektnummer 274762653. 

\bibliographystyle{abbrv}

\end{document}